\documentclass[12pt,draft]{article}

\usepackage{amsfonts,amsmath,amssymb,theorem}

\usepackage[english]{babel}

\DeclareMathOperator{\Id}{Id}
\DeclareMathOperator{\ap}{ap}
\DeclareMathOperator{\adj}{adj}
\DeclareMathOperator{\cp}{cap}
\DeclareMathOperator{\Lip}{Lip}
\DeclareMathOperator{\ACL}{ACL}

\DeclareMathOperator{\loc}{loc}
\DeclareMathOperator{\const}{const}

\DeclareMathOperator{\dist}{dist}
\DeclareMathOperator{\diam}{diam} 
\DeclareMathOperator{\esssup}{esssup}
\DeclareMathOperator{\ess}{ess}
 
\DeclareMathOperator{\BVL}{BVL} 

\begin{document}

\centerline{\bf Sobolev homeomorphisms and Poincar\'e inequality}
\vskip 0.5cm

\centerline{\bf V.~Gol'dshtein and A.~Ukhlov}
\vskip 0.5cm

\centerline{ABSTRACT} 
\bigskip
{\small We study global regularity properties of  Sobolev homeomorphisms on $n$-dimensional Riemannian manifolds under the assumption of  $p$-integrability of its first weak derivatives in degree $p\geq n-1$. We prove that inverse homeomorphisms have integrable first weak derivatives. For the case $p>n$ we obtain necessary conditions for existence of Sobolev homeomorphisms between manifolds. These necessary conditions based on Poincar\'e type inequality:
$$
\inf_{c\in \mathbb R} \|u-c\mid L_{\infty}(M)\|\leq K \|u\mid L^1_{\infty}(M)\|.
$$
As a corollary we obtain the following geometrical necessary condition: 

{\em If there exists a Sobolev homeomorphisms $\varphi: M \rightarrow M'$, $\varphi\in W^1_p(M, M')$, $p>n$, $J(x,\varphi)\ne 0$ a.~e. in $M$, of compact smooth Riemannian manifold $M$ onto Riemannian manifold $M'$ then the manifold $M'$ has finite geodesic diameter.} } 
\vskip 0.5cm

\centerline{\bf Introduction}

\bigskip

Let $D$ and $D'$ be open subsets of Euclidean space $\mathbb R^n$. The classic inverse function theorem implies: {\it Suppose that a homeomorphism $\varphi : D\to D'$ belongs to class $C^k(D)$, $k\geq 1$, and for any point $x \in D$ its differential $D\varphi (x):\mathbb R^n \rightarrow \mathbb R^n$ is an isomorphism. Then the inverse mapping $\varphi^{-1} : D'\to D$ belongs to the same class $C^k(D')$.} 
Here $C^k(D)$ is the vector space of $k$-times continuously differentiable mappings.
 
If a homeomorphism $\varphi$ does not have classic derivatives at all points of $D$, the problem becomes more difficult. In this case the problem can be reformulated in the terms of weakly differentiable mappings and Sobolev spaces.

Let us recall, that Sobolev space
$W^1_p(D)$, $1\leq p<\infty$, consists of locally summable, weakly differentiable functions $f:D\to\mathbb R$ with the finite norm:
$$
\|f\mid W^1_p(D)\|=
\|f\mid\ L_p(D)\|+ \| |\nabla f|\mid L_p(D)\|,\quad \nabla f = (\frac{\partial f}{\partial x_1},...,\frac{\partial f}{\partial x_n}).
$$

Lebesgue space $L_p(D)$, $1\leq p<\infty$, is the space of locally summable functions with the finite norm:
$$
\| f\mid L_p(D)\|=\biggr(\int\limits_D |f|^p~dx\biggl)^{\frac{1}{p}},\quad 1\leq p<\infty.
$$

A mapping $\varphi:D\to\mathbb R^n$ belongs to class $W^1_p(D)$, $1\leq p<\infty$, if
its coordinate functions $\varphi_j$ belong to $W^1_p(D)$, $j=1,\dots,n$.
In this case formal Jacobi matrix
$D\varphi(x)=\bigl(\frac{\partial \varphi_i}{\partial x_j}(x)\bigr)$, $i,j=1,\dots,n$,
and its Jacobian determinant $J(x,\varphi)=\det D\varphi(x)$ are well defined at
almost all points $x\in D$. The norm $|D\varphi(x)|$ of the matrix
$D\varphi(x)$ is the norm of the corresponding linear operator $D\varphi (x):\mathbb R^n \rightarrow \mathbb R^n$ defined by the matrix $D\varphi(x)$. We will use the same notation for this matrix and for the corresponding linear operator.

Necessity of studying of differential and topological properties of mappings inverse to mappings of Sobolev classes arises in problems of the non-linear elasticity theory [1,\,2].
In these works J.~M.~Ball introduced classes of mappings, defined on bounded domains $D \in \mathbb R^n$:
$$
A^+_{p,q}(D)=\{\varphi\in W^1_p(D) : \adj D\varphi\in L_q(D),\quad J(x,\varphi)>0 \quad\text{a.~e. in}\quad D\},
$$
$p,q>n$, where $\adj D\varphi$ is the formal adjoint matrix to the Jacobi matrix $D\varphi$:
$$
\adj D\varphi(x)\cdot D\varphi(x) = \Id J(x,\varphi).
$$
Here $\Id$ is the identical matrix.

The global invertibility was proved for $A^+_{p,q}$-mappings. Topological and weak differential properties of  the inverse mappings were studied.
Fine topological and weak differential properties of these classes of invertible mappings were studied under assumptions 
$p>n-1$, $q\geq {p}/({p-1})$ by Sverak [3] and under less restrictible assumptions $p>n-1$, $q\geq {n}/({n-1})$ By Muller, Tang and Yan [4].

In the frameworks of the inverse function theorem for weakly differentiable mappings the following problem seems to be natural: {\it let $\varphi : D\to D'$ be a Sobolev homeomorphism, what is regularity properties of the inverse homeomorphism.} 

The regularity problem for invertible Sobolev mappings is closely connected with the generalized quasiconformal mappings. Let $D$ and $D'$ be domains in Euclidean space $\mathbb
R^n$, $n\geq 2$. A mapping $\varphi: D\to D'$ is called a {\it mapping
with bounded $(p,q)$-distortion}, $1\leq q\leq p\leq\infty$, if $\varphi$ belongs to Sobolev
space $W^1_{1,\loc}(D)$ and the local $p$-distortion
$$
K_p(x)=\inf\{ k: |D\varphi|(x)\leq k |J(x,\varphi)|^{\frac{1}{p}}, \,\,x\in D\}
$$
belongs to Lebesgue space $L_{r}(D)$, where ${1}/{r}={1}/{q}-{1}/{p}$ (if $p=q$ then $r=\infty$).

Mappings with bounded $(p,q)$-distortion have
a finite distortion, i.~e.  $D\varphi(x)=0$ for almost all points $x$
that belongs to set $Z=\{x\in D:J(x,\varphi)=0\}$.

The class of mappings with bounded $(p,q)$-distortion is a natural generalization of mappings with bounded distortion and represents a non-homeomorphic case of so-called $(p,q)$-quasi\-con\-for\-mal mappings introduced in [5,\,6]. Such classes of mappings have applications to Sobolev type embedding problems [7--9].

\vskip 0.3 cm

The main result of the first section is the next assertion about global regularity of the inverse mapping.

\vskip 0.3 cm

{\bf Theorem A.}
Let a homeomorphism $\varphi : D\to D^{\prime}$, $\varphi(D)=D^{\prime}$ of the Sobolev class
$W^1_{n-1}(D)$ has Luzin $N$-property (the image of a set measure zero is a set measure zero) and $J(x,\varphi)\ne 0$ a.~e. in $D$. Then the inverse
homeomorphism $\varphi^{-1}:D'\to D$  belongs to the Sobolev space $L^1_1(D')$, and for every function $f\in L^1_{\infty}(D)$ the composition $f\circ\varphi^{-1}$ belongs to $L^1_1(D')$.
\vskip 0.3 cm

The next example shows sharpness of the global integrability condition for the inverse homeomorphism.

\vskip 0.3 cm

{\bf Example A.}
For any natural number $k$ consider a family of functions $y=x^{2k+1}$, $-1<x<1$. Their derivatives $y'=(2k+1)x^{2k}$ belong to Lebesgue space $L_p(-1,1)$ for every $p\geq 1$. Inverse functions are $x=y^{\frac{1}{2k+1}}$ and their derivatives $x'=\frac{1}{2k+1}y^{-\frac{2k}{2k+1}}$ belong to $L_1(-1,1)$ and does not belongs to $L_p(-1,1)$ for any fixed apriory $p>1$ and for sufficiently large natural number $k$.
\vskip 0.3 cm

The regularity problem of inverse homeomorphisms to Sobolev
homeomorphisms was studied by many authors. In article [4]
was proved: {\em if a mapping $\varphi\in W^1_{n,\loc}(D)$ and
$J(x,\varphi)>0$ for almost all points $x\in D$, then $\varphi^{-1}$ belongs
to $W^1_{1,\loc}(D')$}.

 \medskip
 
{\em The assumption that $J(x,\varphi)\ne 0$ a.~e. in $D$ cannot be dropped out}. Indeed, consider the function $g(x)=x+u(x)$
on the real line, where $u$ is the standard Cantor function. Let
$f=g^{-1}$. Then the derivative $f'=0$ on the set of positive measure and $h^{-1}$  fails to be absolutely continuous. In this case we can prove only that the inverse
homeomorphism has a finite variation on almost all lines [10]. 
In work [10] was obtained the following result: {\em if a
homeomorphism $\varphi:D\to D'$ belongs to the Sobolev space
$L^1_p(D)$, $p>n-1$, then the inverse mapping $\varphi^{-1}: D'\to D$ has a finite variation on almost all lines (belongs to $\BVL(D')$)}.

\medskip

In work [11] the local regularity of plane homeomorphisms
that belong to Sobolev space $W^1_1(D)$ was studied. For the case of space $\mathbb R^n$, $n\geq 3$, recent work [12] contains the following result for domains in $\mathbb R^n$, $n\geq 3$:
{\em if the norm of the derivative $|D\varphi|$ belongs to Lorentz
space $L^{n-1,1}(D)$ and a mapping $\varphi: D\to D'$ has finite
distortion, then the inverse mapping belongs to Sobolev space
$W^1_{1,\loc}(D')$ and has finite distortion.} Recall that
$$
L^{n-1}(D)\subset L^{n-1,1}(D)\subset \bigcap\limits_{p>n-1}L^p(D).
$$

Note, that results about \emph {local regularity} of inverse mappings for $p=n-1$ were obtained independently in [13] without the assumption about $N$-property of Sobolev homeomorphisms.

\medskip

In the second part of this article we study the problem of equivalence of Riemannian manifolds under Sobolev homeomorphisms of class $W^1_p$. This problem  goes back to well-known Riemann Mapping Theorem, which states existence of a conformal mapping of any plane simply connected domain $D\subset \mathbb R^2$ with non-empty boundary  onto the unit disc $B\subset \mathbb R^2$. From the other side, for Euclidean spaces $\mathbb R^n$, $n\geq 3$, by Liouville theorem the class of conformal mappings is poor and is equal to M\"obius group. Hence, in a natural way arises a problem about description of $n$-dimensional Riemannian manifolds, which are images of  smooth compact $n$-dimensional Riemannian manifolds under homeomorphisms of Sobolev class $W^1_p$. 
A suitable language for this problem is Poncar\'e type inequalities on Riemannian manifolds and its geometrical interpretation.
The main result of this part of the paper is the following:

\vskip 0.3 cm
{\bf Theorem B.}
Let $M$ be a smooth compact $n$-dimensional Riemannian manifold, $M'$ be a smooth $n$-dimensional Riemannian manifold. If a homeomorphism $\varphi : M\to M'$ belongs to $L^1_p(M, M')$, $p>n$, $J(x,\varphi)\ne 0$ a.~e. in $M$, then the manifold $M'$ has a finite geodesic diameter.
\vskip 0.3 cm

Recall, that a homeomorphism $\varphi : D\to D'$ is called quasiconformal, if  local $n$-distortion $K_n$ is uniformly bounded almost everywhere in $D$. It is known, that first weak derivatives of the quasiconformal mapping are locally integrable in degree $n+\varepsilon$, where $\varepsilon$ depends on the distortion $K_n$ [14, 15]. The problem of the global integrability of first weak derivatives of quasiconformal mappings was studied in [16].

We obtain the following assertion about globally integrable in degree $p>2$  conformal mappings.

\vskip 0.3 cm
{\bf Theorem C.}
For any $\varepsilon>0$ does not exist a conformal mapping $\varphi\in L^1_{2+\varepsilon}(B)$ of the unit disc $B\subset\mathbb R^2$ onto a bounded simply connected domain $D\subset\mathbb R^2$ with infinite geodesic diameter.
\vskip 0.3 cm

The similar result is correct for quasiconformal mappings: {\it for any $\varepsilon>0$ does not exist a quasiconformal mapping $\varphi\in L^1_{n+\varepsilon}(B)$ of the unit ball $B\subset\mathbb R^n$ onto a bounded domain $D\subset\mathbb R^n$ with infinite geodesic diameter.}

The suggested method of investigation is based on a relation between Sobolev mappings,  composition operators of spaces of Lipschitz functions and a change of variable formula for weakly differentiable mappings.

\bigskip

\centerline{\bf 1.~Regularity of Sobolev homeomorphisms}

\bigskip

A locally integrable function $f:D\to \mathbb R$ is {\it absolutely continuous on a
straight line $l$} having non-empty intersection with $D$ if it
is absolutely continuous on an arbitrary segment of this
 line which is contained in $D$. A function $f:D\to \mathbb
R$ belongs to the class $\ACL(D)$ ({\it absolutely continuous on
almost all straight lines}) if it is absolutely continuous on
almost all straight lines parallel to any coordinate axis.

Note that $f$ belongs to Sobolev space $L^1_{1}(D)$ if
and only if $f$ is locally integrable and it can be changed by a standard procedure on a set of measure zero (changed to its Lebesgue values at any point where the Lebesgue values exist) so , that a modified
function belongs to
$\ACL(D)$, and its partial derivatives $\frac{\partial f}{\partial
x_i}(x)$, $i=1,\ldots,n$, existing almost everywhere, are
integrable in $D$. From this point we will use modified functions only.
Note that first weak derivatives of
the function $f$ coincide almost everywhere with the usual
partial derivatives (see, e.g., [17] ).

A mapping $\varphi:D\to \mathbb R^n$ belongs to the class $\ACL(D)$, if
its coordinate functions $\varphi_j$ belong to $\ACL(D)$, $j=1,\dots,n$.

We will use the notion of approximate differentiability.
Let $A$ be a subset of $\mathbb R^n$. Density of set $A$ at a point $x\in \mathbb R^n$ is the limit
$$
\lim\limits_{r\to 0}\frac{|B(x,r)\cap A|}{|B(x,r)|}.
$$
Here by symbol $|A|$ we denote Lebesgue measure of the set $A$.

A linear mapping $L: \mathbb R^n\to \mathbb R^n$ is called an approximate differential of a mapping $\varphi : D\to \mathbb R^n$
at point $a\in D$, if for every $\varepsilon>0$ the density of the set
$$
A_{\varepsilon}=\{x\in D : |\varphi(x)-\varphi(a)-L(x-a)|<\varepsilon |x-a| \}
$$
at point $a$ is equal to one.

A point $y\in \mathbb R^n$ is called an approximate limit of a mapping $\varphi : D\to \mathbb R^n$ at a point $x$, if the density of the set $D\setminus\varphi^{-1}(W)$ at this point is equal to zero for every neighborhood $W$ of the point $y$. 

For a mapping $\varphi : D\to\mathbb R^n$ we define approximate partial derivatives
$$
\ap\frac{\partial \varphi_i}{\partial x_j}(x)=\ap\lim\limits_{t \to 0}\frac{\varphi_i(x+te_j)-\varphi_i(x)}{t},\quad i,j=1,...,n.
$$

Approximate differentiable mappings are closely connected with Lipschitz mappings.
Recall, that a mapping $\varphi : D\to \mathbb R^n$ is a Lipschitz mapping, if there exists a constant $K<+\infty$ such that
$$
|\varphi(x)-\varphi(y)|\leq K |x-y|
$$
for every points $x,y \in D$.

The value 
$$
\|\varphi\mid\Lip (D)\| =\sup\limits_{x,y\in D} \frac{|\varphi(x)-\varphi(y)|}{|x-y|}
$$
we call the norm of $\varphi$ in the space $\Lip (D)$.

The next assertion describes this connection between approximate differentiable mappings and Lipschitz mappings in details [18].

\vskip 0.3cm
{\bf Theorem~1.}
Let $\varphi : D\to\mathbb R^n$ be a measurable mapping. Then the following assertions are equivalent:

\noindent
1) The mapping $\varphi : D\to\mathbb R^n$ is approximate differentiable almost everywhere in $D$.

\noindent
2) The mapping $\varphi : D\to\mathbb R^n$ has  approximate partial derivatives $\ap\frac{\partial \varphi_i}{\partial x_j}$, $i,j=1,...,n$ almost everywhere in $D$.

\noindent
3) There exists a collection of closed sets $\{A_k\}_{k=1}^{\infty}$, $A_k\subset A_{k+1}\subset D$, such that a restriction $\varphi \vert_{A_k}$ is a Lipschitz mapping on the set $A_k$ and 
$$
\biggl|D\setminus\sum\limits_{k=1}^{\infty}A_k\biggr|=0.
$$

\vskip 0.3cm

If a mapping $\varphi:D\to D^{\prime}$ has approximate partial
derivatives $\ap\frac{\partial \varphi_i}{\partial x_j}$ almost everywhere
in $D$, $i,j=1,\dots,n$, then  the formal Jacobi matrix
$D\varphi(x)=(\ap\frac{\partial \varphi_i}{\partial x_j}(x))$, $i,j=1,\dots,n$,
and its Jacobian determinant $J(x,\varphi)=\det D\varphi(x)$ are well defined at
almost all points of $D$. The norm $|D\varphi(x)|$ of the matrix
$D\varphi(x)$ is the norm of the linear operator determined by the
matrix in Euclidean space $\mathbb R^n$.

In the theory of mappings with bounded mean distortion additive set functions play a significant role.
Let us recall that a nonnegative mapping $\Phi$ defined on open subsets
of $D$ is called a {\it finitely
quasiadditive} set function [19] if

1) for any point $x\in D$, there exists $\delta$,
$0<\delta<\dist(x,\partial D)$, such that $0\leq
\Phi(B(x,\delta))<\infty$ (here and in what follows
$B(x,\delta)=\{y\in\mathbb R^n: |y-x|<\delta\}$);

2) for any  finite collection $U_i\subset U\subset D$,
$i=1,\dots,k$ of mutually disjoint open sets the following
inequality $\sum\limits_{i=1}^k \Phi(U_i)\leq \Phi(U)$ takes
place.

Obviously, the last inequality can be extended to a countable collection of mutually
disjoint open sets from $D$, so a finitely quasiadditive set
function is also {\it countable quasiadditive.}

If instead of the second condition we suppose that for any finite
collection $U_i\subset D$, $i=1,\dots,k$ of mutually disjoint
open subsets of $D$ the equality
$$
\sum\limits_{i=1}^k \Phi(U_i)= \Phi(U)
$$
takes place, then such set function is said to be {\it finitely
additive}. If the last equality can be extended to a
countable collection of mutually disjoint open subsets of $D$,
then such set function is said to be {\it countable additive.}

A nonnegative mapping $\Phi$ defined on open subsets of $D$ is called a {\it monotone} set function [19] if
$\Phi(U_1)\leq\Phi(U_2)$ under the condition, that $U_1\subset
U_2\subset D$ are open sets.

Note, that a monotone (countable) additive set function is the (countable) quasiadditive set function.

Let us reformulate an auxiliary result from [19] in a convenient for this study way.

\vskip 0.3cm

{\bf Proposition~1.} Let a monotone finitely additive set
function $\Phi$ be defined on open subsets of the domain
$D\subset\mathbb R^n$. Then for almost all points $x\in D$ the
volume derivative
$$
\Phi'(x)=\lim\limits_{\delta\to 0, B_{\delta}\ni x}
\frac{\Phi(B_{\delta})}{|B_{\delta}|}
$$
is finite and for any open set $U\subset D$, the inequality
$$
\int\limits_{U}\Phi'(x)~dx\leq \Phi(U)
$$
is valid.
\vskip 0.3cm

A nonnegative finite valued set function $\Phi$ defined on a collection of
measurable subsets of an open set $D$ is
said to be {\it absolutely continuous} if for every number
$\varepsilon>0$ can be found a number $\delta>0$ such that
$\Phi(A)<\varepsilon$ for any measurable sets $A\subset D$ from
the domain of definition of $\Phi$, which satisfies the condition
$|A|<\delta$.

Let $E$ be a measurable subset of $\mathbb R^n$, $n\geq 2$. Define Lebesgue space $L_p(E)$, $1\leq p\leq\infty$, as a Banach space of locally summable  functions $f:E\to \mathbb R$ equipped with the following norm:
$$
\|f\mid L_p(E)\|=
\biggr(\int\limits_E|f|^p(x)\,dx\biggr)^{1/p},\,\,\,1\leq p<\infty,
$$
and
$$
\|f\mid L_{\infty}(E)\|=
\ess\sup\limits_{x\in E}|f(x)|,\quad p=\infty.
$$

A function $f$ belongs to the space $L_{p,\loc}(E)$, $1\leq p\leq\infty$, if $f\in L_{p}(F)$ for every compact set $F\subset E$.

For an open subset $D\subset\mathbb R^n$ define the seminormed Sobolev space $L^1_p(D)$, $1\leq p\leq\infty$,
as a  space  of locally summable, weakly differentiable functions $f:D\to\mathbb R$ equipped with the following seminorm:
$$
\|f\mid L^1_{p}(D)\|=\| \nabla f\mid L_p(D)\|, \,\,\,\,1\leq p<\infty,  
$$
and
$$
\|f\mid L^1_{\infty}(D)\|=\limsup_{y\to x,\,\,x\in D }\frac{|f(x)-f(y)|}{|x-y|}.  
$$
Here $\nabla f$ is the weak gradient of the function $f$, i.~e. $ \nabla f = (\frac{\partial f}{\partial x_1},...,\frac{\partial f}{\partial x_n})$,
$$
\int\limits_D f \frac{\partial \eta}{\partial x_i}~dx=-\int\limits_D \frac{\partial f}{\partial x_i} \eta~dx, \quad \forall \eta\in C_0^{\infty}(D),\quad i=1,...,n.
$$ 
As usual $C_0^{\infty}(D)$ is the space of infinitely smooth functions with a compact support.

The Sobolev space $W^1_p(D)$, $1\leq p\leq\infty$, is a Banach space of locally summable, weakly differentiable functions $f:D\to\mathbb R$,  equipped with the following norm:
$$
\|f\mid W^1_p(D)\|= \|f\mid L_p(D)\|+\|f\mid L^1_p(D)\|.
$$

A function $f$ belongs to the space $L^1_{p,\loc}(D)$ ($W^1_{p,\loc}(D)$), $1\leq p\leq\infty$, if $f\in L^1_{p}(K)$ ($f\in W^1_{p}(K)$) for every compact subset $K\subset D$.

The Sobolev space $\overset\circ{L}^1_p(D)$ is the
closure of the space $C^{\infty}_0(D)$ in ${L}^1_p(D)$.

A mapping $\varphi:D\to D^{\prime}$ belongs to Lebesgue class $L_p(E)$ if
its coordinate functions $\varphi_j$, $j=1,\dots,n$ belong to $L_p(E)$.
A mapping $\varphi:D\to D^{\prime}$ belongs to Sobolev class $W^1_p(D)$ ($L^1_p(D)$) if
its coordinate functions $\varphi_j$, $j=1,\dots,n$, belong to $W^1_p(D)$ ($L^1_p(D)$).

A mapping $\varphi : D\to D'$ generates by the composition rule $\varphi^{\ast}f=f\circ\varphi$ a bounded composition operator 
$$
\varphi^{\ast} : L^1_p(D')\to L^1_q(D), \,\,\,1\leq q\leq p\leq\infty,
$$
if for every function $f\in L^1_p(D')$ the composition $f\circ\varphi \in L^1_q(D)$ and the inequality
$$
\|\varphi^{\ast}f \mid L^1_q(D)\|\leq K \|f\mid L^1_p(D')\|
$$
holds.
\vskip 0.3cm
{\bf Theorem~2.} Let a homeomorphism $\varphi : D\to D^{\prime}$
between two domains $D,D^{\prime} \subset \mathbb R^n $ belongs to Sobolev space $L^1_{p}(D)$, $1\leq p<+\infty$. Then there exists a bounded
monotone countable additive function $\Phi(A')$ defined on open bounded
subsets of $D'$ such that for every function $f\in
\overset{\circ}{L}^1_{\infty}(A')$ the inequality
$$
\int\limits_{\varphi^{-1}(A)}|\nabla(f\circ \varphi)|^p~dx\leq
\Phi(A'){\esssup\limits_{y\in A'}|\nabla f|^p(y)}
$$
holds.
\vskip 0.3cm

{\sc Proof.} Let us define $\Phi(A')$ by the following way [6]
$$
\Phi(A^{\prime})=\sup\limits_{f\in
\overset{\circ}{L}_{\infty}^{1}(A^{\prime})} \Biggl(
\frac{\bigl\|\varphi^{\ast} f\mid {L}_{p}^{1}(D)\bigr\|}
{\bigl\|f\mid
\overset{\circ}{L}_{\infty}^{1}(A^{\prime})\bigr\|}
\Biggr)^{p},
$$

At first we prove that $\Phi$ is a bounded function. Let a function $f\in L^1_{\infty}(D')$. Then
\begin{multline}
\|\varphi^{\ast}f\mid L^1_p(D)\|=\biggl(\int\limits_D|\nabla (f\circ\varphi)|^p~dx\biggr)^{\frac{1}{p}}
\leq\biggl(\int\limits_D |D\varphi)|^p|\nabla f|^p(\varphi(x))~dx\biggr)^{\frac{1}{p}}\\
\leq\biggl(\int\limits_D |D\varphi)|^p~dx\biggr)^{\frac{1}{p}}\|f\mid L^1_{\infty}(D')\|=
\|\varphi\mid L^1_p(D)\|\cdot \|f\mid L^1_{\infty}(D')\|.
\nonumber
\end{multline}
Hence, we obtain boundedness of the set function $\Phi$ as a consequence of boundedness of the operator $\varphi^{\ast}$.

Let $A_1^{\prime}\subset A_2^{\prime}$ be bounded open
subsets of $D^{\prime}$. Extending
functions of space
$\overset{\circ}{L}_{\infty}^{1}(A_1^{\prime})$ by zero on the set
$A_2^{\prime}$, we obtain an inclusion 
$\overset{\circ}{L}_{\infty}^{1}(A_1^{\prime})\subset
\overset{\circ}{L}_{\infty}^{1}(A_2^{\prime})$. Obviously 
$$
\|f\mid \overset{\circ}{L}_{\infty}^{1}(A_1^{\prime})\|=\|f\mid \overset{\circ}{L}_{\infty}^{1}(A_2^{\prime})\|
$$
for every $f\in \overset{\circ}{L}_{\infty}^{1}(A_1^{\prime})$. By  the following
inequality
\begin{multline}
\Phi(A_1^{\prime})=\sup\limits_{f\in
\overset{\circ}{L}_{\infty}^{1}(A_1^{\prime})} \Biggl(
\frac{\bigl\|\varphi^{\ast} f\mid {L}_{p}^{1}(D)\bigr\|}
{\bigl\|f\mid
\overset{\circ}{L}_{\infty}^{1}(A_1^{\prime})\bigr\|}
\Biggr)^{p}
=
\sup\limits_{f\in
\overset{\circ}{L}_{\infty}^{1}(A_1^{\prime})} \Biggl(
\frac{\bigl\|\varphi^{\ast} f\mid {L}_{p}^{1}(D)\bigr\|}
{\bigl\|f\mid
\overset{\circ}{L}_{\infty}^{1}(A_2^{\prime})\bigr\|}
\Biggr)^{p}\\
\leq\sup\limits_{f\in
\overset{\circ}{L}_{\infty}^{1}(A_2^{\prime})} \Biggl(
\frac{\bigl\|\varphi^{\ast} f\mid {L}_{p}^{1}(D)\bigr\|}
{\bigl\|f\mid
\overset{\circ}{L}_{\infty}^{1}(A_2^{\prime})\bigr\|}
\Biggr)^{p}=\Phi(A_2^{\prime}).
\nonumber
\end{multline}
 the set function $\Phi$ is monotone.

Let $A_{i}^{\prime}$, $i\in \mathbb N$, be open disjoint subsets at the domain
$D^{\prime}$,
$A^{\prime}_0=\bigcup\limits_{i=1}^{\infty}A_i^{\prime}$.
Choose arbitrary functions $f_i\in\overset{\circ}{L}_{\infty}^{1}(A_i^{\prime})$ with following
properties
$$
\bigl\|\varphi^{\ast} f_i\mid{L}_{p}^{1}
(D)\bigr\|\geq \bigl(\Phi(A_i^{\prime})\bigl(1-\frac{\varepsilon}{2^i}\bigr)\bigr)^{\frac{1}{p}}
\bigl\|f_i\mid\overset{\circ}{L}_{\infty}^{1}(A_i^{\prime})
\bigr\|
$$ 
and
$$
\bigl\|f_i\mid\overset{\circ}{L}_{\infty}^{1}(A_i^{\prime})
\bigr\|=1,
$$
while $i\in\mathbb N$. Here $\varepsilon\in(0,1)$ is a fixed number.
Letting $g_N=\sum\limits_{i=1}^{N}f_i$ we obtain
\begin{multline}
\bigl\|\varphi^{\ast} g_N\mid{L}_{p}^{1}
(D)\bigr\|\geq \biggl(\sum\limits_{i=1}^{N}
\left(\Phi(A_i^{\prime})\left(1-\frac{\varepsilon}{2^i}\right)\right)
\bigl\|f_i\mid\overset{\circ}{L}_{\infty}^{1}(A_i^{\prime})
\bigr\|^p\biggr)^{1/p}\\
= \biggl(\sum\limits_{i=1}^{N}\Phi(A_i^{\prime})
\left(1-\frac{\varepsilon}{2^i}\right)\biggr)^{\frac{1}{p}}
\biggl\|g_N\mid\overset{\circ}{L}_{\infty}^{1}
\Bigl(\bigcup\limits_{i=1}^{N}A_i^{\prime}\Bigr)\biggr\|
\\
\geq \biggl(\sum\limits_{i=1}^{N}\Phi(A_i^{\prime})
-\varepsilon\Phi(A_0^{\prime}) \biggr)^{\frac{1}{p}}
\biggl\|g_N\mid\overset{\circ}{L}_{\infty}^{1}
\Bigl(\bigcup\limits_{i=1}^{N}A_i^{\prime}\Bigr)\biggr\|
\nonumber
\end{multline}
since sets, on which the gradients $\nabla\varphi^{\ast} f_i$
do not vanish, are disjoint. From the last inequality  follows that
$$
\Phi(A_0^{\prime})^{\frac{1}{p}}\geq\sup\frac
{\bigl\|\varphi^{\ast} g_N\mid{L}_{p}^{1} (D)\bigr\|}
{\biggl\|g_N\mid\overset{\circ}{L}_{\infty}^{1}
\Bigl(\bigcup\limits_{i=1}^{N}A_i^{\prime}\Bigr)\biggr\|}\geq
\biggl(\sum\limits_{i=1}^{N}\Phi(A_i^{\prime})-\varepsilon\Phi(A_0^{\prime})
\biggr)^{\frac{1}{p}}.
$$
Here the upper bound is taken over all above-mentioned
functions $$
g_N\in\overset{\circ}{L}_{\infty}^{1}
\Bigl(\bigcup\limits_{i=1}^{N}A_i^{\prime}\Bigr).
$$
Since both
$N$ and $\varepsilon$ are arbitrary, we have finally
$$
\sum\limits_{i=1}^{\infty}\Phi(A^{\prime}_i) \leq
\Phi\Bigl(\bigcup\limits_{i=1}^{\infty}A^{\prime}_i\Bigr).
$$

The validity of the inverse inequality can be proved in a straightforward manner.
Indeed, choose functions $f_i\in\overset{\circ}{L}_{\infty}^{1}(A_i^{\prime})$ such that $\bigl\|f_i\mid\overset{\circ}{L}_{\infty}^{1}(A_i^{\prime})
\bigr\|=1$.

Letting $g=\sum\limits_{i=1}^{\infty}f_i$ we obtain
$$
\bigl\|\varphi^{\ast} g\mid{L}_{p}^{1}(D)\bigr\|\leq \biggl(\sum\limits_{i=1}^{\infty}
\Phi(A_i^{\prime})
\bigl\|f_i\mid\overset{\circ}{L}_{\infty}^{1}(A_i^{\prime})
\bigr\|^p\biggr)^{1/p}
= \biggl(\sum\limits_{i=1}^{\infty}\Phi(A_i^{\prime})
\biggr)^{\frac{1}{p}}\biggl\|g_N\mid\overset{\circ}{L}_{\infty}^{1}
\Bigl(\bigcup\limits_{i=1}^{\infty}A_i^{\prime}\Bigr)\biggr\|,
$$
since sets, on which the gradients $\nabla\varphi^{\ast} f_i$
do not vanish, are disjoint. From this inequality follows that
$$
\Phi\biggl(\bigcup\limits_{i=1}^{\infty}A_i^{\prime}\biggr)^{\frac{1}{p}}\leq\sup\frac
{\bigl\|\varphi^{\ast} g\mid{L}_{p}^{1} (D)\bigr\|}
{\biggl\|g\mid\overset{\circ}{L}_{\infty}^{1}
\Bigl(\bigcup\limits_{i=1}^{\infty}A_i^{\prime}\Bigr)\biggr\|}\leq
\biggl(\sum\limits_{i=1}^{\infty}\Phi(A_i^{\prime})
\biggr)^{\frac{1}{p}},
$$
where the upper bound is taken over all functions $g\in\overset{\circ}{L}_{\infty}^{1}
\Bigl(\bigcup\limits_{i=1}^{\infty}A_i^{\prime}\Bigr)$.

By the definition of the set function $\Phi$ we have
$$
\|\varphi^{\ast} f\mid L^1_p(D)\|^p\leq \Phi(A') \|f\mid\overset{\circ}{L}_{\infty}^{1}(A')\|^p
$$
Since the support of the function $f\circ\varphi$ is contained in the set $\varphi^{-1}(A')$ we have
$$
\int\limits_{\varphi^{-1}(A)}|\nabla(f\circ \varphi)|^p~dx\leq
\Phi(A'){\esssup\limits_{y\in A'}|\nabla f|^p(y)}.
$$

Theorem proved.

\vskip0.3cm

We recall some basic facts about $p$-capacity. Let $G\subset\mathbb
R^n$ be an open set and $E\subset G$ be a compact set. For $1\leq
p\leq\infty$ the $p$-capacity of the ring $(E,G)$ is defined as
$$
\cp_p(E,G)=\inf\bigl\{\int_G |\nabla u|^p : u\in
{L}^1_p(G)\cap C^{\infty}_0(G),\, u\geq 1\,\, \text{on}
\,\,E\bigr\}.
$$
Functions $u\in {L}^1_p(G)\cap C^{\infty}_0(G),\, u\geq 1\,\, \text{on}
\,\,E$, are called admissible functions for ring $(E,G)$.

We need the following estimate of the $p$-capacity [20].

\vskip 0.3cm
{\bf Lemma~1.} Let $E$ be a connected closed subset of
an open bounded set $G\subset\mathbb R^n, n\geq 2$, and $n-1<p<\infty$. Then
$$
\cp_p^{n-1}(E,G)\geq c\frac{(\diam E)^p}{|G|^{p-n+1}},
$$
where a constant $c$ depends on $n$ and $p$ only.

\vskip 0.3cm

For readers convenience we will prove this fact.
 
{\sc Proof.} Let $d$ be diameter of set $E$. Without loss of generality we can suggest, that $d=\dist(0, a)$ for some point $a=(0,..,0, a_n)$. For arbitrary number $t$, $0<t<d$, denote by $P_t$ the hyperplane $x_n=t$.

In the subspace $x_n=0$ we consider the unit $(n-2)$-dimensional sphere $S^{n-2}$ with the center at the origin
and fix an arbitrary point $z\in E\cap P_t$. For every point $y\in S^{n-2}$ denote by $R(y)$ the supremum of numbers $r_0$ such that $z+ry\in G$ while $0\leq r\leq r_0$.
Then for every admissible function $f\in C^{\infty}_0(G)$ the following inequality
$$
1= f(z)-f(z+R(y)y)\leq  \int\limits_0^{R(y)}|\nabla f(z+ry)|~dr=\int\limits_0^{R(y)}(|\nabla f(z+ry)|r^{\frac{n-2}{p}})r^{-\frac{n-2}{p}}~dr
$$
holds.
Applying H\"older inequality to the right side of the last inequality, we have
$$
1\leq \biggl(\frac{p-1}{p-n+1}\biggr)^{p-1}\bigl(R(y)\bigr)^{p-n+1}\int\limits_0^{R(y)}|\nabla f(z+ry)|^p r^{n-2}~dr.
$$

Multiplying both sides of this inequality on 
$((p-1)/(p-n+1))^{1-p}\cdot (R(y))^{n-p-1}$ 
and integrating by $y\in S^{n-2}$, we obtain
\begin{multline}
\biggl(\frac{p-1}{p-n+1}\biggr)^{p-1}\int\limits_{S^{n-2}} \bigl(R(y)\bigr)^{p-n+1}~dy\\
\leq\int\limits_{S^{n-2}}~dy
\int\limits_0^{R(y)}|\nabla f(z+ry)|^p r^{n-2}~dr\leq\int\limits_{P_t}|\nabla f|^p~dz.
\nonumber
\end{multline}

For the lower estimate of the left integral we use again H\"older inequality. Denote by $\omega_{n-2}$ the $n-2$-dimensional area of sphere $S^{n-2}$. By simple calculations we get
\begin{multline}
\omega_{n-2}^p=\biggl(\int\limits_{S^{n-2}}~dy\biggr)^p\leq \biggl(\int\limits_{S^{n-2}}\bigl(R(y)\bigr)^{n-p-1}~dy\biggr)^{n-1}
\biggl(\int\limits_{S^{n-2}}\bigl(R(y)\bigr)^{n-1}~dy\biggr)^{p+1-n}\\
\leq((n-1)m_{n-1}(G\cap P_t))^{p-n+1}\biggl(\int\limits_{S^{n-2}}\bigl(R(y)\bigr)^{n-p-1}~dy\biggr)^{n-1}.
\nonumber
\end{multline}
Here $m_{n-1}(A)$ is $(n-1)$-Lebesgue measure of the set $A$.

Denote by $u(t)=m_{n-1}(G\cap P_t)$. Using the last estimate we obtain
$$
\int\limits_{P_t} |\nabla f|^p~dz\geq \biggr( \frac{p-1}{p-n+1}\biggl)^{1-p}(n-1)^{\frac{n-p-1}{n-1}}\omega_{n-2}^{\frac{p}{n-1}}\bigl(u(t)\bigr)^{\frac{n-p-1}{n-1}}.
$$
After integrating by $t\in (0,d)$ we have
$$
\int\limits_{G} |\nabla f|^p~dx\geq \biggr( \frac{p-1}{p-n+1}\biggl)^{1-p}(n-1)^{\frac{n-p-1}{n-1}}\omega_{n-2}^{\frac{p}{n-1}}\int\limits_0^d\bigl(u(t)\bigr)^{\frac{n-p-1}{n-1}}~dt.
$$
By H\"older inequality
\begin{multline}
d^p=\biggl(\int\limits_0^d~dt\biggr)^p\leq\biggl(\int\limits_0^du(t)~dt\biggr)^{p-n+1} 
\biggl(\int\limits_0^d\bigl(u(t)\bigr)^{\frac{n-p-1}{n-1}}~dt\biggr)^{n-1} \\
\leq |G|^{p-n+1}\biggl(\int\limits_0^d\bigl(u(t)\bigr)^{\frac{n-p-1}{n-1}}~dt\biggr)^{n-1}.
\nonumber
\end{multline}

Therefore
$$
\int\limits_{G}|\nabla f|^p~dx\geq \biggr( \frac{p-1}{p-n+1}\biggl)^{1-p}(n-1)^{\frac{n-p-1}{n-1}}\omega_{n-2}^{\frac{p}{n-1}}\biggl(\frac{d^p}{|G|^{p-n+1}}\biggr)^{\frac{1}{n-1}}.
$$
Since $f$ is an arbitrary admissible function the required inequality is proved.

\vskip 0.3cm

Let us define a class $BVL$ of mappings with finite variation.
A mapping $\varphi : D\to\mathbb R^n$ belongs to the class $\BVL(D)$
(i.e., has {\it finite variation on almost all straight lines})
if it has finite variation on almost all straight lines $l$ parallel
to any coordinate axis: for any finite number of points $t_1,...,t_k$ that belongs to such straight line $l$
$$
\sum\limits_{i=0}^{k-1}|\varphi(t_{i+1})-\varphi(t_i)|<+\infty.
$$

For a mapping $\varphi$ with finite variation on almost all straight lines, the partial
derivatives $\partial \varphi_i/\partial x_j$, $i,j=1,\dots,n$,
exists almost everywhere in $D$.

\vskip 0.3cm
{\bf Theorem~3. [10]} Let a homeomorphism $\varphi : D\to D^{\prime}$
belongs to Sobolev space $L^1_{p}(D)$, $p>n-1$. Then the inverse
homeomorphism $\varphi^{-1}:D'\to D$ belongs to the class $\BVL(D')$.
\vskip 0.3cm

For readers convenience we reproduce here a slightly modified proof of this fact.

{\sc Proof.} Take an arbitrary $n$-dimensional open
parallelepiped $P$ such that $\overline{P}\subset D'$ and its
edges are parallel to coordinate axis. Let us show that
$\varphi^{-1}$ has finite variation on almost all intersection of $P$
and straight lines parallel to $x_n$-axis.

Let $P_0$ be the projection of $P$ on the subspace $x_n=0$, and
let $I$ be the projection of $P$ on the coordinate axis $x_n$.
Then $P=P_0\times I$. The monotone countable-additive function
$\Phi$ determines a monotone countable additive function of open
sets $A\subset P_0$ by the rule $\Phi(A,P_0)=\Phi(A\times I)$. For
almost all points $z\in P_0$, the quantity
$$
\overline{\Phi'}(z,P_0) =\overline{\lim_{r\to
0}}\biggl[\frac{\Phi(B^{n-1}(z,r),P_0)} {r^{n-1}}\biggr]
$$
is finite [21] (here $B^{n-1}(z,r)$ is the $(n-1)$-dimensional ball of
radius $r>0$ centered at the point $z$).

The $n$-dimensional Lebesgue measure $\Psi(U)=|\varphi^{-1}(U)|$, where
$U$ is an open sen in $D'$, is a monotone countable additive
function and, therefore, also determines a monotone countable
additive function $\Psi(A,P_0) = \Psi(A\times I)$ defined on open
sets $A\subset P_0$. Hence  $\overline{\Psi'}(z,P_0)$ is finite for
almost all points $z\in P_0$.

Choose an arbitrary point $z\in P_0$ where 
 $\overline{\Phi'}(z,P_0)<+\infty$ and $\overline{\Psi'}(z,P_0)<+\infty$. On the section $I_z=\{z\}\times I$ of the
parallelepiped $P$, take arbitrary mutually disjoint closed
intervals $\Delta_1,...,\Delta_k$ with lengths $b_1,...,b_k$
respectively. Let $R_i$ denote the open set of points for which
distances from $\Delta_i$ smaller than a given $r>0$:
$$
R_i=\{x\in G : \dist(x,\Delta_i)<r\}.
$$
Consider
the ring $(\Delta_i, R_i)$. Let $r>0$ be selected so that
$r<cb_i$ for $i=1,\dots,k$, where $c$ is a sufficiently small constant.
Then the function $u_i(x)=\dist(x,\Delta_i)/r$ is an admissible
for ring $(\Delta_i, R_i)$. 

By Theorem~2 we have the estimate
$$
\|\varphi^{\ast} u_i \mid L^1_p(D)\|^p\leq \Phi(A') \|u_i \mid\overset{\circ}{L}_{\infty}^{1}(A')\|^p
$$
for every function $u_i$, $i=1,...,k$.

Hence, for every ring $(\Delta_i, R_i)$, $i=1,...,k$, the inequality
$$
\cp_p(\varphi^{-1}(\Delta_i), \varphi^{-1}(R_i))\leq \Phi(R_i)\cp_{\infty}(\Delta_i, R_i)
$$
holds. 

The function $u_i(x)=\dist(x,\Delta_i)/r$ is admissible
for ring $(\Delta_i, R_i)$ and we have the upper estimate
$$
\cp_{\infty}(\Delta_i, R_i)\leq |\nabla u_i|=\frac{1}{r}.
$$

Applying the
lower bound for the capacity of the ring (Lemma~1), we obtain
$$
\biggl(\frac{(\diam
\varphi^{-1}(\Delta_i))^{p/(n-p)}}{|\varphi^{-1}(R_i)|^{(p-n+1)(n-1)}}\biggr)^{\frac{1}{p}}\leq
c_1 \Phi(R_i)^{\frac{1}{p}}\cdot\frac{1}{r}.
$$
This inequality gives
$$
\diam \varphi^{-1}(\Delta_i)\leq
c_2\biggl(\frac{|\varphi^{-1}(R_i)|}{r^{n-1}}\biggr)^{\frac{p-n+1}{p}}
\cdot \biggl(\frac{\Phi(R_i)}{r^{n-1}}\biggr)^{\frac{n-1}{p}}.
$$
Summing over $i=1,\dots,k$
we obtain
$$
\sum\limits_{i=1}^k\diam \varphi^{-1}(\Delta_i)\leq
c_2 \sum\limits_{i=1}^k \biggl(\frac{|\varphi^{-1}(R_i)|}{r^{n-1}}\biggr)^{\frac{p-n+1}{p}}
\cdot \biggl(\frac{\Phi(R_i)}{r^{n-1}}\biggr)^{\frac{n-1}{p}}.
$$

Hence
$$
\sum\limits_{i=1}^k\diam \varphi^{-1}(\Delta_i)\leq
c_2 \biggl(\sum\limits_{i=1}^k\frac{|\varphi^{-1}(R_i)|}{r^{n-1}}\biggr)^{\frac{p-n+1}{p}}
\cdot \biggl(\sum\limits_{i=1}^k\frac{\Phi(R_i)}{r^{n-1}}\biggr)^{\frac{n-1}{p}}.
$$

Using the Besicovitch type theorem [22] for the estimate of the
value of the function $\Phi$ in terms of the multiplicity of a
cover, we obtain
$$
\sum\limits_{i=1}^k\diam \varphi^{-1}(\Delta_i)\leq
c_3 \biggl(\frac{|\varphi^{-1}(\bigcup_{i-1}^k R_i)|}{r^{n-1}}\biggr)^{\frac{p-n+1}{p}}
\cdot \biggl(\frac{\Phi(\bigcup_{i-1}^k R_i)}{r^{n-1}}\biggr)^{\frac{n-1}{p}}.
$$
Hence
$$
\sum\limits_{i=1}^k \diam \varphi^{-1}(\Delta_i)\leq
c_3\biggl(\frac{|\varphi^{-1}(B^{n-1}(z,r),P_0)|}{r^{n-1}}\biggr)^{\frac{p-n+1}{p}}
\cdot
\biggl(\frac{\Phi(B^{n-1}(z,r),P_0)}{r^{n-1}}\biggr)^{\frac{n-1}{p}}.
$$
Because  $\overline{\Phi'}(z,P_0)<+\infty$ and $\overline{\Psi'}(z,P_0)<+\infty$ we obtain finally 

$$
\sum\limits_{i=1}^k \diam \varphi^{-1}(\Delta_i)<+\infty.
$$
Therefore $\varphi^{-1}\in \BVL(D')$.

Theorem proved.

\vskip 0.3cm

Let us recall the change of variable formula for Lebesgue integral [23].
Let a mapping $\varphi : D\to \mathbb R^n$ be such that
there exists a collection of closed sets $\{A_k\}_1^{\infty}$, $A_k\subset A_{k+1}\subset D$ for which restrictions $\varphi \vert_{A_k}$ are Lipschitz mapping on sets $A_k$ and 
$$
\biggl|D\setminus\sum\limits_{k=1}^{\infty}A_k\biggr|=0.
$$
Then there exists a measurable set $S\subset D$, $|S|=0$ such that  the mapping $\varphi:D\setminus S \to \mathbb R^n$ has Luzin $N$-property and the change of variable formula
$$
\int\limits_E f\circ\varphi (x) |J(x,\varphi)|~dx=\int\limits_{\mathbb R^n\setminus \varphi(S)} f(y)N_f(E,y)~dy
$$ 
holds for every measurable set $E\subset D$ and every nonnegative Borel measurable function $f: \mathbb R^n\to\mathbb R$. Here 
$N_f(y,E)$ is the multiplicity function defined as the number of preimages of $y$ under $f$ in $E$.

If a mapping $\varphi$ possesses Luzin $N$-property (the image of a set of measure zero has measure zero), then $|\varphi (S)|=0$ and the second integral can be rewritten as the integral on $\mathbb R^n$.
Note, that if a homeomorphism $\varphi : D\to D'$ belongs to the Sobolev space $W^1_{n,\loc}(D)$ then $\varphi$ has Luzin $N$-property and the change of variable formula holds [24].

If a mapping $\varphi : D\to \mathbb R^n$ belongs to the Sobolev space $W^1_{1,\loc}(D)$ then by [23] there exists a collection of closed sets $\{A_k\}_1^{\infty}$, $A_k\subset A_{k+1}\subset D$ for which restrictions $f\vert_{A_k}$ are Lipschitz mapping on sets $A_k$ and 
$$
\biggl|D\setminus\sum\limits_{k=1}^{\infty}A_k\biggr|=0.
$$ 
Hence for such mappings the previous change of variable formula is correct.

Like in [25] (see also [12]) we define a measurable function

$$
\mu(y)=
\begin{cases}
\biggl(\frac{|\adj D\varphi|(x)}{|J(x,\varphi)|}\biggr)_{x=\varphi^{-1}(y)}
\quad\text{if}\quad J(x,\varphi)\neq 0,
\\
\,\,0\quad\quad\quad\quad\quad\quad \quad\quad \quad    \text{otherwise}.
\end{cases}
$$

Because $J(x,\varphi)\neq 0$ almost everywhere and $\varphi$ has Luzin $N$-property, the function $\mu(y)>0$ almost everywhere in $D'$.

The following lemma was proved but does not formulated in [12] under an additional assumption that $|D\varphi|$ belongs to the Lorentz space $L^{n-1,n}(D)$.

\vskip 0.3cm
{\bf Lemma~2.} Let a homeomorphism $\varphi:D\to D', \varphi(D)=D'$ belongs to the
Sobolev space $L^1_p(D)$ for some $p\geq n-1$. Then the function $\mu$ is locally
integrable in the domain $D'$.
\vskip 0.3cm

{\sc Proof.} 
Using the change of variable formula for Lebesgue integral [23] and Luzin $N$-property of $\varphi$  we have the following equality
$$
\int\limits_{D'}\mu(y)~dy=\int\limits_{D'\setminus \varphi(S)}\mu(y)~dy=\int\limits_{D\setminus S} | \mu(\varphi(x))|J(x,\varphi)|~dx= \int\limits_{D} |\adj D\varphi|(x)~dx.
$$
Applying H\"older inequality, we obtain that for every compact subset $F'\subset D'$
$$
\int\limits_{F'}\mu(y)~dy\leq \int\limits_{F} |\adj D\varphi|(x)~dx\leq C \int\limits_{F}
|D\varphi|^{n-1}(x)~dx,
$$
where $F'=\varphi(F)$.
Therefore,   $\mu$ belongs to $L_{1,\loc}(D')$, since $\varphi$ belongs to $L^1_p(D)$, $p\geq n-1$, and as consequence $\varphi\in L^1_{n-1,\loc}(D)$.
\vskip 0.3cm

{\bf Theorem~4.} Suppose that a homeomorphism $\varphi : D\to D^{\prime}$, $\varphi(D)=D^{\prime}$ belongs to Sobolev space
$W^1_{p}(D)$, $p\geq n-1$, has Luzin $N$-property (the image of a set measure zero is a set measure zero) and $J(x,\varphi)\ne 0$ a.~e. in $D$. Then the inverse
homeomorphism $\varphi^{-1}:D'\to D$  has integrable first weak derivatives, is a
mapping with bounded $(p',1)$-distortion for $p'=p/(p-n+1)$ and induces a bounded composition operator
$$
(\varphi^{-1})^{\ast} : L^1_{p'}(D)\to L^1_1(D'),\quad (\varphi^{-1})^{\ast}(f)=f\circ\varphi^{-1}.
$$
\vskip 0.3cm

{\bf Corollary~1.} Suppose a homeomorphism $\varphi : D\to D^{\prime}$ 
belongs to Sobolev space $L^1_{n-1}(D)$, has Luzin $N$-property and $J(x,\varphi)\ne 0$ a.~e. in $D$. Then the inverse mapping $\varphi^{-1}:D'\to D$ belongs to the Sobolev space $L^1_1(D')$
\vskip 0.3cm

{\sc Proof of Corollary~1.}
By Theorem~4  the operator
$$
(\varphi^{-1})^{\ast} : L^1_{\infty}(D)\to L^1_1(D'),\quad (\varphi^{-1})^{\ast}(f)=f\circ\varphi^{-1}
$$
is bounded.

Substituting in the inequality
$$
\|f\circ \varphi^{-1}\mid L^1_1(D')\|\leq K \|f\mid L^1_{\infty}(D)\|
$$
coordinate functions $x_j$, $j=1,...,n$, we conclude, that coordinate functions $\varphi^{-1}_j$ belongs to $L^1_1(D')$. Hence, $\varphi^{-1}\in L^1_1(D')$.
\vskip 0.3cm

{\sc Proof of Theorem~4.} Since absolute continuity is the local property, it is sufficient to prove that the mapping $\varphi^{-1}$ belongs to $\ACL$ on every compact subset of $D'$. Consider arbitrary cube $Q'\in D'$, $\overline{Q'}\in D$, with edges parallel to coordinate axes, and $Q=\varphi^{-1}(Q')$. For $i=1,\dots n$ we will use a notation: $Y_i=(x_1,...,x_{i-1},x_{i+1},...,x_n)$,
$$
F_i(x)=(\varphi_1(x),\dots,\varphi_{i-1}(x),\varphi_{i+1}(x),\dots,\varphi_n(x))
$$
and  $Q'_i$  is the intersection of the cube $Q'$ with a line
$Y_i=\const$. 

Using the change of variable formula and the Fubini theorem [26] we obtain the following estimate
$$
\int\limits_{F_i(Q)}H^{n-1}(dY_i)\int\limits_{Q_i'}\mu(y)~H^1(dy)=\int\limits_{Q'} \mu(y)~dy=\int\limits_Q |\adj D\varphi| (x)~dx<+\infty.
$$
Hence for almost all $Y_i\in F_i(Q)$
$$
\int\limits_{Q_i'}\mu(y)~H^1(dy)<+\infty.
$$

Let $\ap J\varphi(x)$ be an approximate Jacobian of
the trace of the mapping $\varphi$ on the set $\varphi^{-1}(Q'_i)$ [26].  Consider a point $x\in Q$ in which there exists a non-generated approximate differential $\ap Df(x)$ of the mapping $\varphi: D\to D'$. Let $L: \mathbb R^n\to \mathbb R^n$ be a linear mapping induced by this approximate differential $\ap Df(x)$. We denote by the symbol $P$ the image of the unit cube $Q_0$ under the linear mapping $L$ and by $P_i$ the intersection of $P$ with the image of the line $x_i=0$. Let $d_i$ be a length of $P_i$. Then 
$$
d_i\cdot|\adj DF_i|(x)=|Q_0|=|J(x,\varphi)|.
$$
So, since $d_i=\ap J\varphi(x)$ we obtain that for
almost all $x\in Q\setminus Z$, $Z=\{x\in D : J(x,\varphi)=0\}$, we have
$$
\ap J\varphi(x)=\frac{|J(x,\varphi)|}{|\adj DF_i|(x)}.
$$

Let $\tilde{D}$ be a subset of $D$ such that the differential $D\varphi$ is well defined at all points of the set $\tilde{D}$. Since $|D\setminus \tilde{D}|=0$, we can assume, without loss of generality, that $D\setminus \tilde{D}=S$, where $S$ is the set from the change of variable formula. By the Luzin $N$-property $|\varphi(S)|=0$.

Now, since $|\varphi(S)|=0$ and the function $\mu$ is equal zero on the set $\varphi(S)$ we have for arbitrary
compact set $A'\subset Q'_i$, and for almost all $Y_i\subset F_i(Q)$,
the following inequality
$$
\int\limits_{A'\setminus \varphi(S)}\mu(y)~H^1(dy)=\int\limits_{A'\setminus \varphi(S)}\mu(y)~H^1(dy)+\int\limits_{A'\cap \varphi(S)}\mu(y)~H^1(dy)=\int\limits_{A'}\mu(y)~H^1(dy)<+\infty.
$$
By construction of $\mu$:
\begin{multline}
 H^1(\varphi^{-1}(A'))\leq \int\limits_{\varphi^{-1}(A')}\frac{|\adj
D\varphi|(x)}{|\adj
DF_i|(x)}~H^1(dx)\\
=\int\limits_{\varphi^{-1}(A')}\frac{|\adj
D\varphi|(x)}{|J(x,\varphi)|}\cdot\frac{|J(x,\varphi)|}{|\adj
DF_i|(x)}~H^1(dx)
=
\int\limits_{\varphi^{-1}(A')\setminus S}\mu(\varphi(x))\cdot\frac{|J(x,\varphi)|}{|\adj
DF_i|(x)}~H^1(dx)\\
+
\int\limits_{\varphi^{-1}(A')\cap S}\mu(\varphi(x))\cdot\frac{|J(x,\varphi)|}{|\adj
DF_i|(x)}~H^1(dx)
=
\int\limits_{\varphi^{-1}(A')}\mu(\varphi(x))\ap
J\varphi(x)~H^1(dx).
\nonumber
\end{multline}
The last equality holds because $\mu(\varphi(x))=0$ on the set $S$.

By using the change of variable formula for the Lebesgue integral [26, 27] we obtain
$$
H^1(f^{-1}(A'))\leq \int\limits_{A'}\mu(y)~H^1(dy).
$$
Therefore, the mapping $\varphi^{-1}$ is absolutely continuous on almost all
lines in $D'$and is a weakly differentiable mapping.

Since the homeomorphism $\varphi$ has Luzin $N$-property then preimage of a set positive measure is a set positive measure. Hence, the volume derivative of the inverse mapping
$$
J_{\varphi^{-1}}(y)=\lim\limits_{r\to 0}\frac{|\varphi^{-1}(B(y,r))|}{|B(y,r)|}>0
$$
almost everywhere in $D'$.
So $J(y,\varphi^{-1}) \neq 0$ for almost all points $y\in D$. Integrability of the $p'$-distortion follows from the inequality
$$
|D\varphi^{-1}|(y)\leq |D\varphi(x)|^{n-1}\big/ |J(x,\varphi)|
$$
which holds for almost all points $y=\varphi(x)\in D'$.

Indeed, with the help of the change of variable formula, we have
\begin{multline}
\int\limits_{D'}\biggl(\frac{|D\varphi^{-1}(y)|^{p'}}{|J(y,\varphi^{-1})|}\biggr)^{\frac{1}{p'-1}}~dy=
\int\limits_{D'}\biggl(\frac{|D\varphi^{-1}(y)|}{|J(y,\varphi^{-1})|}\biggr)^{\frac{p'}{p'-1}}|J(y,\varphi^{-1})|~dy\\
\leq \int\limits_{D}\biggl(\frac{|D\varphi^{-1}(\varphi(x))|}{|J(\varphi(x),\varphi^{-1})|}\biggr)^{\frac{p'}{p'-1}}~dx\leq \int\limits_{D}|D\varphi(x)|^p~dx<+\infty.
\nonumber
\end{multline}

The boundedness of the composition operator follows from integrability of the $p'$-distortion [6].

If $p=n-1$ we have the following estimate
$$
\int\limits_{D'}|D\varphi^{-1}(y)|~dy\leq
\int\limits_{D}|D\varphi^{-1}(\varphi(x))||J(x,\varphi)|~dx
\leq \int\limits_{D}|D\varphi(x)|^{n-1}~dx<+\infty.
\nonumber
$$

Hence $\varphi^{-1}\in L^1_1(D')$.

Theorem proved.
\vskip 0.3cm

For higher integrability of the inverse mapping we need some additional assumptions.
Let us recall
that a homeomorphism $\varphi$ between domains $D$ and $D'$ is said to
be $(p,q)$-quasiconformal mapping [6] if one of the following
equivalent assertions is fulfilled:

1) the mapping $\varphi$ belongs to $W^1_{1,\loc}(D)$ and the value
$$
K_{p,q}(D)=\|K_p(\cdot)\mid L_{\kappa}(D)\|
$$
is finite, where $K_p(x)=\inf\{ k: |D\varphi|(x)\leq k |J(x,\varphi)|^{\frac{1}{p}},\,\, x\in D\}$. Here number $\kappa$ is defined from the
relation $1/\kappa = 1/q-1/p$.

2) the mapping $\varphi$ induces a bounded embedding operator of
Sobolev spaces
$$
\varphi^{\ast} : L^1_p(D')\to L^1_q(D), \quad 1\leq q\leq p<\infty,
$$
by the rule $\varphi^{\ast}(u)=u\circ \varphi$. Moreover, the norm of the
operator is equivalent to the value $K_{p,q}(D)$.

\vskip 0.3cm

{\bf Theorem~5.} Let $\varphi : D\to D^{\prime}$ be a $(p,q)$-quasiconformal homeomorphism, $n-1\leq q\leq p\leq \infty$, $\varphi$ has Luzin $N$-property and $J(x,\varphi)\ne 0$ a.~e. in $D$ while $q=n-1$. Then the inverse mapping $\varphi^{-1} : D'\to D$ is the $(q',p')$-quasiconformal homeomorphism, $p'=p/(p-n+1)$, $q'=q/(q-n+1)$. 
\vskip 0.3cm

{\sc Proof.} In the case $\infty>p\geq q>n-1$ the theorem was proved in [6]. Let $q=n-1$, then $q'=\infty$. Because $J(x,\varphi)\ne 0$ a.~e. in $D$ and homeomorphism $\varphi$ has Luzin $N$-property we can conclude that $J(y,\varphi^{-1}) \neq 0$ for almost all $y\in D'$. Then by the change of variable formula
\begin{multline}
\int\limits_{D'} |D\varphi^{-1}|^{p'}(y)~dy=\int\limits_{D'\setminus\varphi(S)} |D\varphi^{-1}|^{p'}(y)~dy 
=\int\limits_{D} |D\varphi^{-1}|^{p'}(\varphi(x))|J(x,\varphi)|~dx\\
\leq \int\limits_{D} \biggl(\frac{|D\varphi|^{n-1}(x)}{|J(x,\varphi)|}\biggr)^{p'}|J(x,\varphi)|~dx=
\int\limits_{D} \biggl(\frac{|D\varphi|^{p}(x)}{|J(x,\varphi)|}\biggr)^{\frac{n-1}{p-(n-1)}}~dx<+\infty,
\nonumber
\end{multline}
since $\varphi$ is the $(p,n-1)$-quasiconformal homeomorphism and has Luzin $N$-property. Hence $\varphi^{-1}: D'\to D$ belongs to the Sobolev space $L^1_{p'}(D)$ and generates the bounded composition operator
$$
\varphi^{\ast} : L^1_{\infty}(D)\to L^1_{p'}(D').
$$
Therefore the inverse mapping $\varphi^{-1}$ is a $(\infty, p')$-quasiconformal homeomorphism.

For $p=\infty>q$ all conditions of Theorem~4 are fulfilled: the homeomorphism $\varphi$ belongs to $L^1_q(D)$, has Luzin $N$-property and $J(x,\varphi)\ne 0$ a.~e. in $D$. Therefore by Theorem 4 the inverse mapping $\varphi^{-1} : D'\to D$ is the $(q',1)$-quasiconformal homeomorphism.
\vskip 0.3cm

If $p=q=\infty$, then the mapping $\varphi$ generates the bounded composition operator
$$
\varphi^{\ast} : L^1_{\infty}(D')\to L^1_{\infty}(D).
$$
Hence, substituting coordinate functions $f_i=y_i$, $i=1,...,n$, we obtain that the mapping $\varphi\in L^1_{\infty}(D)$ and is differentiable almost everywhere in $D$. Therefore at almost all points $y=\varphi(x)\in D'$ and at almost all points $x\in D\setminus Z$, $Z=\{x\in D : J(x,\varphi)=0\}$,
we have 
$$
|D\varphi^{-1}|(y)\leq \frac{|D\varphi|^{n-1}(x)}{|J(x,\varphi)|}={|D\varphi|^{n-1}(x)}\cdot |J(y,\varphi^{-1})|
\leq \|\varphi \mid L^1_{\infty}(D)\|\cdot |J(y,\varphi^{-1})|,
$$
and $\varphi^{-1}$ is $(1,1)$-quasiconformal homeomorphism.
\vskip 0.3cm
From Theorem~5 immediately follows that under conditions of this theorem for $p>n-1$
the following fact is correct.

\vskip 0.3cm

{\bf Corollary~2.} Let a homeomorphism $\varphi : D\to D^{\prime}$, $\varphi\in W^1_{n-1,\loc}(D)$, has Luzin $N$-property, $J(x,\varphi)\ne 0$ a.~e. in $D$, and
$$
\int\limits_{D} \biggl(\frac{|D\varphi|^{p}(x)}{|J(x,\varphi)|}\biggr)^{\frac{n-1}{p-(n-1)}}~dx<+\infty.
$$
Then the inverse mapping $\varphi^{-1} : D'\to D$ belongs to $L^1_{p'}(D')$. 
\vskip 0.3cm

\bigskip

\centerline{\bf 2. About existence of Sobolev mappings}

\bigskip

Let $M$ be a $n$-dimensional Riemannian manifold and $\omega$ be a differential $k$-form defined on $M$. 

Recall the notion of the weak exterior differential of the differential form $\omega$. We denote by $C_0^{\infty}(\mathbb M, \Lambda^k)$ the vector space of smooth differential forms of degree $k$ with compact support on $\mathbb M$ and by $L_{\loc}^{1}(\mathbb M, \Lambda^k)$ the space of differential $k$-forms whose coefficients (in any local coordinate system) are locally integrable. We say that a form $\theta\in L_{\loc}^{1}(\mathbb M, \Lambda^k)$ is the weak exterior differential of the differential form $\omega\in L_{\loc}^{1}(\mathbb M, \Lambda^{k-1})$ if for each $\phi\in C_0^{\infty}(\mathbb M, \Lambda^{n-k})$
$$
\int\limits_{\mathbb M}\theta\wedge \phi=(-1)^k\int\limits_{\mathbb M}\omega\wedge  d\phi.
$$
We denote the weak exterior differential of the form $\omega$ by $d\omega$. Remember that the weak exterior differential of a function $f$ is a 1-form $df$. This notion of weakly differentiable functions coincides with 
the standard definition for Euclidean domains.  

Define Sobolev space $L^1_p(M)$, $1\leq p\leq\infty$, as a space of locally integrable functions $f : M\to\mathbb R$ which have a weak exterior differential $df$ and the finite seminorm:
$$
\|f\mid L^1_p(M)\| = \||df|\mid L_p(M)\|.
$$

Let $M$ and $M'$ be $n$-dimensional Riemannian manifolds. We say, that a mapping $\varphi : M\to M'$ belongs to the Sobolev class $L^1_p(M, M')$, if for every Lipschitz function $f : M'\to \mathbb R$ the composition $f\circ\varphi$ belongs to the Sobolev class $L^1_p(M)$. 

The norm of the differential $|D\varphi|$ is defined:
$$
|D\varphi|(x)=\sup\limits_{f\in \Lip(M')}\frac{|d(f\circ\varphi)|}{\|f\mid\Lip(M')\|}
$$

Note, that in the case $M=M'=\mathbb R^n$ this definition is equivalent to the standard definition of the weak differential for  mappings of Euclidean domains. Therefore results of the previous section are correct for the Sobolev homeomorphisms defined on Riemannian manifolds.

In this section we study existence conditions for Sobolev homeomorphisms of a smooth compact $n$-dimensional Riemannian manifold $M$ onto a Riemannian manifold $M'$ in terms of Poincar\'e type inequalities existence on $M'$.
For this purpose we will use the following "`transfer"' diagram for Poincar\'e inequalities (see, [7--9]):

\begin{eqnarray}
L^1_{\infty}({M}')\overset{\varphi^{\ast}}{\longrightarrow}L^{1}_{p}(M)
\nonumber
\\
\hskip -6cm\downarrow \hskip 2.5cm \downarrow
\nonumber
\\
{ L_{\infty}({M}')\overset{(\varphi^{-1})^{\ast}} {\longleftarrow}L_{\infty}(M)}
\nonumber
\end{eqnarray}
Here the operator $\varphi^{\ast}$ defined by the composition rule $\varphi^{\ast}f=f\circ\varphi$ is a bounded composition operator of Sobolev spaces induced by a homeomorphism $\varphi$ of manifolds $M$ and $M'$.

\vskip 0.3cm

{\bf Theorem~6.} 
Let $M$ be a smooth compact Riemannian $n$-dimensional manifold, $M'$ be a Riemannian $n$-dimensional manifold. Suppose that a homeomorphism $\varphi : M\to M'$ belongs to $L^1_p(M, M')$, for some $p>n$, and $J(x,\varphi)\ne 0$ a.~e. in $M$. Then every function $f\in L^1_{\infty}(M')$  is a continuous function on $M'$ and the following inequality
\begin{equation}
\inf_{c\in \mathbb R} \|f-c\mid L_{\infty}(M')\|\leq K \|f\mid L^1_{\infty}(M')\|
\end{equation}
holds.

Here constant $K>0$ depends only on a mapping $\varphi$ and a manifold $M$.
\vskip 0.3cm

{\sc Proof.} Choose arbitrarily $f \in L^1_{\infty}(M')$. By the definition of the class $L^1_p(M, M')$ a composition $f\circ\varphi$ belongs to Sobolev space $L^1_p(M)$. Since $M$ is a smooth compact Riemannian $n$-dimensional manifold, then for the space $L^1_p(M)$ the Poincar\'e type inequality
$$
\inf_{c\in \mathbb R} \|f\circ\varphi-c\mid L_{\infty}(M)\|\leq K' \|f\circ\varphi\mid L^1_{p}(M)\|
$$
is fullfiled and by the classical embedding theorem $f\circ\varphi \in C(M)$.

Since $\varphi^{-1}$ is a homeomorphism then $f=(f\circ\varphi)\circ\varphi^{-1} \in C(M')$  and
$$
\|f-c\mid L_{\infty}(M')\|=\|f\circ\varphi-c\mid L_{\infty}(M)\|
$$
for every $c\in\mathbb R$.

From the other side, since $J(x,\varphi)\ne 0$ a.~e. in $M$, we have that
\begin{multline}
\|f\circ\varphi\mid L^1_{p}(M)\|=\biggl(\int\limits_M |\nabla f\circ\varphi|^p~dx\biggr)^{\frac{1}{p}}\leq
\biggl(\int\limits_M |\nabla f|^p(\varphi(x)) |D\varphi|^p~dx\biggr)^{\frac{1}{p}}\\
\leq \| |\nabla f| (\varphi(x)) \mid L_{\infty}(M)\|\biggl(\int\limits_M |D\varphi|^p~dx\biggr)^{\frac{1}{p}}
=\| |\nabla f| \mid L_{\infty}(M')\|\biggl(\int\limits_M |D\varphi|^p~dx\biggr)^{\frac{1}{p}}.
\nonumber
\end{multline}

Therefore
$$
\inf_{c\in \mathbb R} \|f-c\mid L_{\infty}(M')\|\leq K' \|\varphi\mid L^1_p(M)\| \|f\mid L^1_{\infty}(M')\|.
$$

Theorem proved.

\vskip 0.3cm

{\bf Example~1.} 
Does not exists a homeomorphism $\varphi \in L^1_p(B, D_1)$, $J(x,\varphi)\ne 0$ a.~e. in $B$, of the unit ball $B \in \mathbb R^n$ onto the domain
$$
D_1=\{(x,y)\in \mathbb R^2 : 1< x<\infty,0< y< 1\}
$$
 for any  $p>n$.
\vskip 0.3cm

{\sc Proof.}
Consider the function $f(x,y)=\frac{1}{\sqrt{2}}(x+y)$, defined on the domain $D_1$. Then $|\nabla f|=1$ and 
$\|\nabla f\mid L_{\infty}({D}_1)\|=1$. From the other side, $f(x,y)\rightarrow \infty$, as $x\rightarrow\infty$. Hence for the function $f(x,y)$ the inequality (1)  doesn't holds and by Theorem 6 for any $p>n$ a homeomorphism $\varphi \in L^1_p(B, D_1)$ can not exists.
\vskip 0.3cm

In the case $p=n$ we have following simple necessary condition for existence of a homeomorphisms $\varphi$ of the class $ W^1_n$.

\vskip 0.3cm

{\bf Lemma~3.} 
Let $M$ be a smooth compact Riemannian $n$-dimensional manifold, $M'$ be a Riemannian $n$-dimensional manifold, and a homeomorphism $\varphi : M\to M'$ belongs to $L^1_n(M, M')$. Then the manifold $M'$ has finite measure.
\vskip 0.3cm

{\sc Proof.} By change of variable formula for Sobolev mappings (see, for example,[28]) and H\"older inequality we have
$$
|M'|=\int\limits_M |J(x,\varphi)|~dx\leq \int\limits_M\biggl(\sum\limits_{i,j=1}^n\bigl(\frac{\partial\varphi_j}{\partial x_i}\bigr)^2\biggr)^{\frac{1}{2}}dx \leq C\|\varphi\mid L^1_n(M, M')\|.
$$
\vskip 0.3cm

In the next example we construct a finite measure domain $D_{\alpha}$ which doesn't allows a homeomorphism of  class $L^1_p(B,D_{\alpha})$, $p>n$ of the unit ball $B\in R^n$ onto $D_{\alpha}$. 

{\bf Example~2.}
Does not exists a homeomorphism $\varphi \in L^1_p(B, D_{\alpha})$, $J(x,\varphi)\ne 0$ a.~e. in $B$, of the unit ball $B \in \mathbb R^n$ onto the domain
$$
{D}_{\alpha}=\{(x,y)\in \mathbb R^2 : 1< x<\infty,0< y< 1/x^{\alpha},\,\,\alpha>1\}.
$$
for any $p>n$.
\vskip 0.3cm

{\sc Proof.}
Let us remark that $|D_{\alpha}|<\infty$ because $\alpha>1$ and 
$$
|D_{\alpha}|=\int\limits_{D_{\alpha}}dxdy=\int\limits_1^{\infty}\frac{dx}{x^{\alpha}}<\infty.
$$

Consider the function $f(x,y)=\frac{1}{\sqrt{2}}(x+y)$, defined on the domain ${D}_{\alpha}$. Then $|\nabla f|=1$ and 
$\|\nabla f\mid L_{\infty}({D}_{\alpha})\|=1$. From the other side, $f(x,y)\rightarrow \infty$, as $x\rightarrow\infty$. Hence for the function $f(x,y)$ the inequality (1)  doesn't holds and by Theorem 6 for any $p>n$ a homeomorphism $\varphi \in L^1_p(B, D_{\alpha})$ can not exists.
\vskip 0.3cm

\vskip 0.3cm

{\bf Theorem~7.} 
Let $M$ be a smooth compact Riemannian $n$-dimensional manifold, $M'$ be a Riemannian $n$-dimensional manifold, and a homeomorphism $\varphi : M\to M'$ belongs to $L^1_p(M, M')$, $p>n$, $J(x,\varphi)\ne 0$ a.~e. in $M$. Then the manifold $M'$ has a finite geodesic diameter.
\vskip 0.3cm

{\sc Proof.}
Let $x_0$ be an arbitrary fixed point at the manifold $M'$. Consider the function $f(x)=\dist_{M'}(x_0,x)$. Here $\dist_{M'}(x,y)$ is the intrinsic geodesic distance.
Then $|\nabla f|=1$ almost everywhere on $M'$ and $f$ belongs to $L^1_{\infty}(M')$. By Theorem~6 the function $f$ belongs to the Lebesgue space $L_{\infty}(M')$ and we have the following estimate:
$$
\ess\sup\limits_{x\in M'}\dist_{M'}(x_0,x)=\|f\mid L_{\infty}(M')\|<+\infty.
$$
\vskip 0.3cm

{\bf Example~3.} {\em Bounded domain with finite measure and with infinite geodesic diameter.}

Consider the plane domain $R=\{(x,y)\in \mathbb R^2 : 1<x^2+y^2<4 \}$. Fix a number $0<r<1/4$ and for natural number $n$, $n=1,2,...$, consider circles:
$$
S^r_{2n}=\biggl\{(x,y)\in \mathbb R^2 : x^2+y^2=\bigl(1+\frac{1}{2n}\bigr)^2, \quad x<1-r \biggr\}
$$
and
$$
S^r_{2n+1}=\biggl\{(x,y)\in \mathbb R^2 : x^2+y^2=\bigl(1+\frac{1}{2n+1}\bigr)^2, \quad x>r-1 \biggr\}.
$$

Let 
$$
D=R\setminus (\cup_{n=1}^{\infty} S^r_{2n})\cup (\cup_{n=1}^{\infty} S^r_{2n+1}).
$$
Then  for every $\varepsilon>0$ there no exists a conformal homeomorphism $\varphi$ of unit disc $B$ onto domain $D$ of the class $L^1_{2+\varepsilon}(B)$. 
\vskip 0.3cm

{\sc Proof.} It is clear that this domain $D$ has a infinite geodesic diameter and finite measure.

We obtain the following assertion about globally integrable in degree $p>2$  conformal mappings as a direct corollary of Theorem 7.

\vskip 0.3 cm
{\bf Theorem C.}
For any $\varepsilon>0$ does not exist a conformal mapping $\varphi\in L^1_{2+\varepsilon}(B)$ of the unit disc $B\subset\mathbb R^2$ onto a bounded simply connected domain $D\subset\mathbb R^2$ with infinite geodesic diameter.
\vskip 0.3 cm

{\bf Remark~1.} Remind that by Riemann mapping Theorem there exists a conformal homeomorphism $\phi : B\to D$ but this conformal homeomorphism can not belong to the class $L^1_{2+\varepsilon}(B)$ for any $\varepsilon>0$. 

\vskip 0.3cm
The similar result is correct for quasiconformal mappings: {\it for any $\varepsilon>0$ does not exist a quasiconformal mapping $\varphi\in L^1_{n+\varepsilon}(B)$ of the unit ball $B\subset\mathbb R^n$ onto a bounded domain $D\subset\mathbb R^n$ with infinite geodesic diameter.}

\vskip 0.5cm
We are pleasure to thank Professor Jan Maly for helpful discussions.

\vskip 0.3cm

\centerline{REFERENCES}

\begin{enumerate}

\item
Ball J.~M. {\it Convexity conditions and existence theorems in nonlinear elasticity.}//
Arch. Rat. Mech. Anal. -- 1976. -- V.~63. -- P.~337--403.

\item
Ball J.~M. {\it Global invertability of Sobolev functions and the interpretation of matter.}//
Proc. Roy. Soc. Edinburgh -- 1981. -- V.~88A. --  P.~315--328.

\item
\~Sver\'ak V. {\it Regularity properties of deformations with finite energy.}//
Arch. Rat. Mech. Anal. -- 1988. -- V.~100. -- P.~105--127.

\item
Muller S., Tang Q., Yan B.~S. {\it On a new class of elastic deformations not allowing for cavitation.}//
Ann. Inst. H. Poincare. Anal. non. lineaire -- 1994. -- V.~11. -- N.~2. -- P.~217--243.

\item
Ukhlov A.~D. {\it Mappings that generate embeddings of Sobolev spaces.}//
Siberain Math. J. -- 1993. -- V.~34. -- N.~1. -- P.~165--171.

\item
Vodop'yanov S.~K., Ukhlov A.~D. {\it Sobolev spaces and $(p,q)$-quasiconformal mappings of Carnot groups.}//
Siberain Math. J. -- 1998. -- V.~39. -- N.~4. -- P.~776--795.

\item
Gol'dshtein V., Gurov L. {\it Applications of change of variable operators for exact embedding theorems.}//
Integral equations operator theory -- 1994. -- V.~19. -- N.~1. -- P.~1--24.

\item
Gol'dshtein V., Ramm A.~G. {\it Compactness of the embedding operators for rough domains.}//
Math. Inequalities and Applications -- 2001. -- V.~4. -- N.~1. -- P.~127--141.

\item
Gol'dshtein V., Ukhlov A. {\it Weighted Sobolev spaces and embedding theorems.}//
Transactions of Amer. Math. Soc. (to appear)

\item
Ukhlov~A. {\it Differential and geometrical properties of Sobolev mappings.}//
Mathematical Notes -- 2004. -- V.~75. -- N.~2. -- P.~291--294.

\item
Hencl S., Koskela P. {\it Regularity of the inverse of a planar Sobolev homeomorphism.}//
Arch. Rational Mech. Anal. -- 2006. -- V.~180. -- N.~1. -- P.~75--95.

\item
Hencl S., Koskela P., Maly Y. {\it Regularity of the inverse of a Sobolev homeomorphism in space.}//
Proc. Roy. Soc. Edinburgh Sect. A. -- 2006. -- V.~136A. -- N.~6. -- P.~1267--1285.

\item
Cs\"ornyei M., Hencl S., Maly Y. {\it Homeomorphisms in the Sobolev space $W^{1,~n-1}$.}//
Preprint. -- MATH-KMA-2007/252.

\item
Boyarski\u\i ~ B.~V. {\it Homeomorphic solutions of Beltrami systems.}//
Dokl. Akad. Nauk SSSR -- 1955. -- V.~102. -- P.~661--664.

\item
Gehring F.~W. {\it The $L^{p}$-integrability of the partial derivatives of quasiconformal mapping.}//
Bull. Amer. Math. Soc. -- 1973. -- V.~79. -- P.~465--466.

\item
Astala K., Koskela P. {\it Quasiconformal mappings and global integrability of the derivative.}//
J. Anal. Math.  -- 1991. -- V.~57. -- P.~203--220.

\item
Maz'ya V. {\it Sobolev spaces} -- Berlin: Springer Verlag. 1985.

\item
Whitney H. {\it On total differentiable and smooth functions.}//
Pacific J. Math. -- 1951. -- N.~1. -- P.~143--159.

\item
Vodop'yanov~S.~K., Ukhlov~A.~D. {\it Set functions and its applications in the theory of Lebesgue and Sobolev spaces.}//
Siberian Adv. Math. -- 2004. -- V.~14. -- N.~4. -- P.~1--48.

\item
Kruglikov V.~I. {\it Capacities of condensers and spatial mappings quasiconformal in the mean.}//
Matem. sborn. -- 1986. -- V.~130. -- N.~2. -- P.~185--206.

\item
Rado T.,Reichelderfer P.~V. {\it Continuous Transformations in Analysis} -- Berlin: Sp\-rin\-ger Verlag. 1955.

\item
Gusman M. {\it Differentiation of integrals in $\mathbb R^n$} -- Moscow: Mir. 1978.

\item
Hajlasz P. {\it Change of variable formula under minimal assumptions.}//
Colloq. Math. -- 1993. -- V.~64. -- N.~1. -- P.~93--101.

\item
Reshetnyak Yu.~G. {\it Some geometrical properties of functions and mappings with generalized derivatives.}//
Siberain Math. J. -- 1966. -- V.~7. -- P.~886--919.

\item
Peshkichev Yu.~A. {\it Inverse mappings for homeomorphisms of the class $BL$.}//
Mathematical Notes -- 1993. -- V.~53. -- N.~5. -- P.~98--101.

\item
Federer H. {\it Geometric measure theory} -- Berlin: Sp\-rin\-ger Verlag. 1969.

\item
Hajlasz P. {\it Sobolev mappings, co-area formula and related topics.}//
Proc. on analysis and geometry. Novosibirsk -- 2000. -- P.~227--254.

\item
Gol'dshtein V.~M., Reshetnyak Yu.~G. {\it Quasiconformal mappings and Sobolev spaces} -- Dordrecht, Boston, London: Kluwer Academic Publishers. 1990.

\end{enumerate}
\end{document}